**The Error Term In the Prime Counting Function**
**N. A. Carella, August, 2009**

*Abstract:* This article considers the error term of the prime counting function. It applies some recent results on the densities of primes in short intervals to derive an improvement of the error term from subexponential size to fractional exponential size. The corresponding equivalent results for the zeros of the zeta function and the Mertens function are also discussed.

**1. Introduction**
For a real number $x > 0$, define the primes counting function by $\pi(x) = \#\{ p \leq x : p \text{ is prime }\}$. This work is an attempt to reduce the error term $E(x)$ in the analytic formula $\pi(x) = li(x) + E(x)$. The oldest formula for the prime counting function appears to be

$$\pi(x) \approx \frac{x}{\log x - A(x)}, \qquad (1)$$

where $\lim A(x) = 1.08366...$, as $x$ tends to infinity. This claim is attributed to Legendre, 1798. Subsequently, the formula

$$\pi(x) \approx \int_2^x \frac{dt}{\log t} \qquad (2)$$

was conjectured by Gauss, 1849, [IG, p. 3], [PT], [BN], et cetera. These approximations are very close and can be written as

$$\pi(x) = \frac{x}{\log(x)} + O\left(\frac{x}{\log^2(x)}\right), \qquad (3)$$

where the error term $E_0(x) = O(x/\log^2 x)$, see [RB, p. 213], [GR]. After foundational works by Chevyschev and Riemann, circa 1850, this claim was independently proved by both DelaValle Poussin and Hadamard, circa 1896. The DelaValle Poussin form of the error term has the shape

$$E_1(x) = \pi(x) - li(x) = O(xe^{-c(\log x)^{1/2}}), \qquad (4)$$

where $li(x)$ denotes the logarithm integral, $c > 0$ is an absolute constant, see [NW] for extensive details. The current version of the error term is works of several authors in the 1950's. The Vinogradov form of the error term has the shape

$$E_2(x) = \pi(x) - li(x) = O(xe^{-c(\log x)^{3/5}(\log \log x)^{-1/5}}), \qquad (5)$$

confer [IV], [FD]. The corresponding form of the prime counting function (Prime Number Theorem) is

$$\pi(x) = li(x) + O(xe^{-c(\log x)^{3/5}(\log \log x)^{-1/5}}). \qquad (6)$$

The optimum error term $E_z(x) = O(x^{1/2} \log x)$ is specified by the theory of the zeta function, and it is also known that $\pi(x) - li(x) = \Omega_{\pm}(x^{1/2}(\log x)^{-1} \log \log \log x)$, Littlewood, 1914.

In addition to the analytical approximations, there are a few exact formulae of the prime counting function. The oldest appears to be the Legendre formula from the 1700's, specifically,

$$\pi(x) = \pi(x^{1/2}) - 1 + \sum_{d \mid p_1 \cdot p_2 \cdots p_n} \mu(d) \left[\frac{x}{d}\right], \tag{7}$$

where $p_n$ is the $n$th prime number and $p_n \leq x^{1/2}$. This expression is derived from the sieve of Eratosthenes, 200 BC. The other prominent exact formula is the vonMangoldt formula

$$\pi(x) = \sum_{n=1}^{\infty} \frac{\mu(n)}{n} \left( li(x^{1/n}) - \sum_{\rho} li(x^{\rho/n}) + \int_{x^{1/n}}^{\infty} (t(t^2-1)\log t)^{-1} dt \right), \tag{8}$$

which is derived from the zeros $\rho$ of the zeta function. This formula was discovered by Riemann in 1850, and proved by vonMangoldt about fifty years later.

An improvement of the current error term of the prime counting function is achieved by mean of recent results [HL], [HB], et cetera, on the densities of prime numbers in short intervals. The new estimate of the error term has fractional exponential size $E(x) = O(x^{7/12-\varepsilon(x)})$ as opposed to subexponential size $E_2(x) = O(xe^{-.218(\log x)^{3/5}(\log\log x)^{-1/5}})$.

**Theorem 1.** Let $x > 1$ be a sufficiently large number. Then

$$\pi(x) = li(x) + O(x^{7/12-\varepsilon(x)}), \tag{9}$$

where $0 \leq \varepsilon(x) < 1/12$ and $\varepsilon(x) \to 0$ as $x \to \infty$.

Assuming the fundamental results on the densities of the primes in short intervals in [HL] and [HB], this is an elementary and easy to verify proof, a proof straight from the 'Book', see [NT].

The error term of the prime counting function, primes in short intervals, the distribution of the zeros of the zeta function on the critical strip, and a few related topics constitute a circle of equivalent mathematical statements. Accordingly, any advance in any of these fundamental topics immediately cascades down to the other equivalent topics. In light of this observation, it is not surprising that the advances on the densities of primes in short intervals in [HL], [HB], [BK], et cetera, lead to equivalent advances on the distribution of the zeros of the zeta function on the critical strip and the error term of the prime counting function and so on. In fact, the result in [BK] can yield a slightly better result than that stated in Theorem 1 above.

The second section introduces the notations and some elementary results. It concludes with the improved versions of the Chebychev theta and psi functions. In the third Section the new improved version of the theta function is utilized to derive the proof of Theorem 1. The fourth Section ends with two applications to the zeros of the zeta function and to Mertens function.

## 2. Elementary Results
Several fundamental results are recorded here for the convenience of the reader. The proofs of these results, often lengthy and difficult, are given in the cited literature.

**2.1 Primes in Short Intervals and Their Densities.** Some relevant claims concerning the prime numbers on short intervals are recalled, further, several fundamental results on the densities of prime numbers on various short intervals are also given.

**Theorem 2.** ([BK]) For all $x > x_0$, the interval $(x - x^{.525}, x]$ contains prime numbers.

The estimate for the maximum number of prime numbers in short intervals is often called Brun-Titchmarsh Theorem. Using state of the art sieve methods, this theorem states the followings.

**Theorem 3.** ([RE]) There exists an $x_0$ such that for all $x \geq x_0$ and all $y \geq 1$,

$$\pi(x+y) - \pi(x) \leq \frac{2y}{\log y + 3.53}. \tag{10}$$

The problem of the local density of the primes (or equivalently of finding a formula for counting the primes in short intervals) and the problem of finding the maximum gap are closely related but not the same. The former is a slightly more complex problem. The simplest case deals with the asymptotic density of prime numbers in the long interval $[2, x]$, this is given by

$$\frac{\pi(x)}{x} \sim \frac{1}{\log(x)}. \tag{11}$$

A challenging aspect of the prime density problem is concerned with the density of primes in the short interval $(x, x+y]$, $0 < y < x$, and their distributions. The basic density formula over the short intervals $(x, x+y]$ should have the asymptotic form

$$\frac{\pi(x+y) - \pi(x)}{y} \sim c_0 \frac{1}{\log(x+y)} \tag{12}$$

and the constant should be $c_0 = 1$, see [M] and [GR] for counterexamples and discussions. It is an important problem to determine the smallest $y > 0$ such that the formula holds.

**Theorem 4.** ([SG]) Let $f(x)$ be positive and increasing and $f(x)/x$ decreasing for $x > 0$, further, suppose that $f(x)/x \to 0$ and $\lim \log f(x)/\log x > 19/77$ for $x \to \infty$. Then for almost all $x > 0$ the density $\frac{\pi(x+f(x)) - \pi(x)}{f(x)} \sim \frac{1}{\log x}$ holds.

The best unconditional density result states the following.

**Theorem 5.** ([HB]) Let $\varepsilon(x) \leq 1/12$ be a nonnegative function of $x$. Then

(i) $\pi(x) - \pi(x-y) = \frac{y}{\log x}\left(1 + O(\varepsilon(x)^4) + O((\log\log x/\log x)^4)\right),$ (13)

(ii) $\pi(x) - \pi(x - x^{7/12}) = \frac{x^{7/12}}{\log x}\left(1 + O((\log\log x/\log x)^4)\right),$

uniformly for $x^{7/12-\varepsilon(x)} \leq y \leq x/(\log x)^4$.

The situation for small $y < x^{1/2}$ is quite different, probably because the density of prime numbers in very short intervals is not as uniform as in larger intervals.

**Theorem 6.** ([M]) Let $f(x) = (\log x)^\delta$, $\delta > 1$. Then

$$\limsup_{x \to \infty} \frac{\pi(x+f(x)) - \pi(x)}{f(x)/\log x} > 1 \quad \text{and} \quad \liminf_{x \to \infty} \frac{\pi(x+f(x)) - \pi(x)}{f(x)/\log x} < 1 \tag{14}$$

infinitely often. For the range $1 < \delta < e^\gamma$, the limit superior satisfies $\limsup_{x \to \infty} \frac{\pi(x+f(x)) - \pi(x)}{f(x)/\log x} \geq \frac{e^\gamma}{\delta}$.

**2.2 The Theta and Psi Functions.** The theta and psi functions are defined by $\vartheta(x) = \sum_{p \leq x} \log p$ and $\psi(x) = \sum_{p^n \leq x, n \geq 1} n \log p$ respectively. Some details on the early history of these estimates appears in [NW]. New results are given in [MT], [DT] etc.

**Theorem 7.** (Chebyshev 1850) Let $x > 1$ be a real number. Then
(i) $ax < \vartheta(x) < bx$ for some constants $a, b > 0$,          (ii) $\psi(x) = \vartheta(x) + O(x^{1/2})$,

Proof: The first can be derived from the properties of the binomial coefficients $(n,k) = n!/(k!(n-k)!)$, see [HW, p. 453], [JS], and the second follows from the finite series

$$\psi(x) = \vartheta(x) + \vartheta(x^{1/2}) + \vartheta(x^{1/3}) + \cdots + \vartheta(x^{1/\log x}) = \vartheta(x) + \vartheta(x^{1/2}) + \vartheta(x^{1/3}) \log x = \vartheta(x) + O(x^{1/2}). \quad \blacksquare$$

The error terms $\vartheta(x) - x$ and $\psi(x) - x$ are oscillating functions that tend to infinity as $x$ tends to infinity. The peaks and valleys of the oscillations of the latter function are known to satisfy $|\psi(x) - x| > c\sqrt{x}$, $c > 0$ constant, infinitely often. Further, the difference $\psi(x) - \vartheta(x) = \vartheta(x^{1/2}) + \cdots + \vartheta(x^{1/n}) = x^{1/2} + O(x^{1/3})$ has a striking resemblance to a *discontinuous* square root function, viz, $\psi(x) - \vartheta(x) \asymp \sqrt{x}$ almost everywhere.

**Theorem 8.** (Littlewood 1914) The error term satisfies $|\psi(x) - x| = \Omega_\pm(x^{1/2} \log \log x)$.

**Theorem 9.** ([MT]) Assuming the Riemann hypothesis, the followings hold.
(i) The function $\psi(x) - x$ changes sign in the interval $[x, 19x]$, $x \geq 1$ at least ounce.
(ii) The function $\psi(x) - x$ changes sign in the interval $[x, 2.02x]$, $x \geq 1$ at least ounce.

Confer [MT] for the analysis. A weaker result on the sign changes is given in [NW, p. 226]. It is believed that

$$\liminf_{x \to \infty} \frac{\psi(x) - x}{x^{1/2}(\log\log\log x)^2} = -\frac{1}{2\pi} \quad \text{and} \quad \limsup_{x \to \infty} \frac{\psi(x) - x}{x^{1/2}(\log\log\log x)^2} = \frac{1}{2\pi}. \tag{15}$$

An explicit formula version of the psi function, which shows the influence of the zeros of the zeta function is included here, see [ES, p. 50] for details.

**Theorem 10.** (Landau Formula) Let $\{ \rho = \sigma + it : 0 \leq \sigma \leq 1, t \in \mathbb{R} \}$ be the set of critical zeros of $\zeta(s)$. Then

$$\psi(x) = x - \sum_{|\rho| \leq T} \frac{x^\rho}{\rho} + O(xT^{-1}(\log Tx)^2) + O(\log x), \tag{16}$$

uniformly for $T > T_0$.

Consider the nontrivial zeros $\sigma + it \in \mathbb{C}$ of the zeta function $\zeta(\sigma + it) = 0$. The maximum $\theta$ of the real part $\sigma$ of the zeros of the zeta function is known to satisfy $0 < \sigma < 1$. Further, by their symmetry about the line $\sigma = 1/2$, it follows that $1/2 \leq \sigma < 1$, see [IG, p. 83].

**Theorem 11.** Let $\theta = \max \{ \sigma : \zeta(\sigma + it) = 0, 0 < \sigma < 1 \}$, and let $x \geq x_0$ be a real number. Then
(i) $\psi(x) = x + O(x^\theta \log^2 x)$,          (ii) $\vartheta(x) = x + O(x^\theta \log^2 x)$,

Proof: Put $T = x$ and use $t_n = an/\log n$, where $\sigma_n + it_n$ is the $n$th zeros, in the Landau formula to show (i), then use Theorem 7 to show (ii), confer [IG, p83]. ∎

The best estimates on the theta and psi functions have subexponentials error terms. These are derived using a zero-free region of the form $\sigma > 1 + c(\log t)^{-2/3}(\log\log t)^{1/3}$, where $c > 0$ is an absolute constant, and the explicit formula, [IK, p. 227].

***Lemma* 12.** Let $x \geq x_0$, then
(i) $\psi(x) = x + O(xe^{-c(\log x)^{3/5}(\log\log x)^{-1/5}})$.
(ii) $\vartheta(x) = x + O(xe^{-c(\log x)^{3/5}(\log\log x)^{-1/5}})$.

Confer [RS], [SC] and [DT] for other estimates. The new estimates of the theta and psi functions spring from the fact that $\vartheta(x) = x + c_0 x / f(x)$ for some function $f(x)$ in tandem with the densities results in [HL] and [HB] for primes in short intervals.

***Theorem* 13.** Let $x > 1$ be a sufficiently large number. Then
(i) $\vartheta(x) = x + ax^{7/12 - \varepsilon(x)}$,
(ii) $\psi(x) = x + bx^{7/12 - \varepsilon(x)}$,

where $a, b > 0$ are constants and $0 \leq \varepsilon(x) < 1/12$ and $\varepsilon(x) \to 0$ as $x \to \infty$.

Proof: Let $x > 1$ be a sufficiently large number and let $y > x^{7/12 - \varepsilon(x)}$. Suppose that $\vartheta(x) = x + c_0 x / f(x)$, where $f(x) = o(x)$ and $c_0 > 0$ is a constant. Three cases depending on the function $f(x)$ will be considered.

Case 1. The constants $c_0 \neq c_1$, and $\vartheta(x) = x + c_0 x / f(x)$, where $f(x) = (\log x)^B$, $B > 0$ is a constant.
Case 2. The constants $c_0 \neq c_1$, and $\vartheta(x) = x + c_0 x / f(x)$, where $f(x) = e^{c(\log x)^\beta}$, $c > 0$ and $0 < \beta < 1$.

These cases lead to contradictions, so $\vartheta(x) \neq x + c_0 x / \log^B(x)$ nor $\vartheta(x) \neq x + c_0 x / e^{c(\log x)^\beta}$. The analysis of all these cases are similar, see the Appendix for cases 1 and 2.

Case 3. The constants $c_0 \neq c_1$, and $\vartheta(x) = x + c_0 x / f(x)$, where $f(x) = x^\tau$, $\tau > 0$ is a small number.
Then the theta difference over the short interval $(x, x + y]$ is given by

$$\vartheta(x+y) - \vartheta(x) = x + y + c_1 \frac{x+y}{(x+y)^\tau} - \left(x + c_0 \frac{x}{x^\tau}\right)$$
$$= y\left(1 + \frac{c_1}{(x+y)^\tau}\right) + x\left(\frac{c_1}{(x+y)^\tau} - \frac{c_0}{x^\tau}\right). \quad (17)$$

Now replace the power series expansion of $(x+y)^\tau = x^\tau(1 + \tau y/x + O(y^2/x^2))$ to simplify the last equation to

$$\vartheta(x+y) - \vartheta(x) = y\left(1 + \frac{c_1}{(x+y)^\tau}\right) + x\left(\frac{c_1 - c_0(1 + \tau y/x + O(y^2/x^2))}{(x+y)^\tau}\right). \quad (18)$$

By Theorem 5, the number of primes in the short interval $(x, x+y]$ with $y \geq x^{7/12 - \varepsilon(x)}$, is within the range

$$\frac{c_2 y}{\log(x+y)} \leq \pi(x+y) - \pi(x) \leq \frac{c_3 y}{\log(x+y)}, \qquad (19)$$

where $c_0, c_1, c_2, c_3, c_4, \ldots$ are constants. These estimates and the identity $\vartheta(x+y) - \vartheta(x) = \sum_{x < p \leq x+y} \log p$ immediately yield

$$\frac{c_2 y}{\log(x+y)} \log x \leq \sum_{x < p \leq x+y} \log p \leq \frac{c_3 y}{\log(x+y)} \log(x+y). \qquad (20)$$

Next replacing (18) into (20), then dividing by $y \geq c_4 x^{7/12 - \varepsilon(x)}$ across the board and simplifying return

$$c_2 \left(1 - \frac{y}{x \log x}\right) \leq 1 + \frac{c_1}{(x+y)^\tau} + \frac{x^{5/12+\varepsilon(x)}}{c_4}\left(\frac{c_1 - c_0(1 + \tau y/x + O(y^2/x^2))}{(x+y)^\tau}\right) \leq c_3. \qquad (21)$$

Since both sides of these inequalities are constants for all sufficiently large number $x > 1$ and a proper choice of constant $c_4 > 0$, it quickly follows that $\tau = 5/12 + \varepsilon(x)$ is the minimum exponent. This proves the claim (i). The claim (ii) follows from (i) and Theorem 7-ii. ∎

***Corollary* 14.** For all sufficiently large numbers $x > 1$, and a random variable $y \in (x^{7/12 - \varepsilon(x)}, x]$, the followings hold.
(i) The random variable $|\vartheta(x+y) - \vartheta(x) - y| \leq ax^{7/12-\varepsilon(x)}$,
(ii) The random variable $|\psi(x+y) - \psi(x) - y| \leq bx^{7/12-\varepsilon(x)}$,

where $a, b > 0$ are constants, and $0 \leq \varepsilon(x) < 1/12$ and $\varepsilon(x) \to 0$ as $x \to \infty$.

The current perspective views the expression $\vartheta(x+y) - \vartheta(x) - y$ as a normal random variable of mean $\mu = 0$ and variance $\sigma = x\log(x/y)$. An abstract of the precise statement as given in [MJ], is as follows.

***Conjecture* 15.** Let $\varepsilon > 0$, $N > 1$, and $N^\varepsilon < y < N^{1-\varepsilon}$. For $y \leq x < N$, the random variable $\psi(x+y) - \psi(x) - y$ is normally distributed with mean $\sim y$ and variance $\int_2^N (\psi(x+y) - \psi(x) - y)^2 dx \sim Ny \log N/y$.

**2.4. Representations of Arithmetic Functions.** This section considers the representations of some arithmetic functions. The representations use integrals and other related functions. These are elementary techniques but very effective in their realm of applications.

***Theorem* 16.** For a real number $x > 1$, the following statements hold.

(i) $\pi(x) = \dfrac{\vartheta(x)}{\log x} + \displaystyle\int_2^x \dfrac{\vartheta(t)}{t \log^2 t} dt$, (ii) $\vartheta(x) = \pi(x) \log x - \displaystyle\int_2^x \dfrac{\pi(t) dt}{t \ln^2 t}$, (22)

(iii) $\displaystyle\sum_{2 \leq n \leq x} \dfrac{\Lambda(n)}{\log n} = \dfrac{\psi(x)}{\log x} + \int_2^x \dfrac{\psi(t) dt}{t \ln^2 t}$, (iv) $\displaystyle\sum_{p \leq x} 1/p = \dfrac{\pi(x)}{x} + \int_2^x \dfrac{\pi(t)}{t^2} dt$.

Proof of (i): There are a few proofs of this identity. One of these proofs proceeds as shown below. Recall that $\vartheta(x) = \sum_{p \leq x} \log p$. This permits to write the line

$$\int_2^x \frac{\vartheta(t)}{t \log^2 t} dt = \int_2^x \frac{\sum_{p \leq t} \log p}{t \log^2 t} dt = \sum_{p \leq x} \log p \int_p^x \frac{1}{t \log^2 t} dt = \sum_{p \leq x} \log p \left(\frac{1}{\log p} - \frac{1}{\log x}\right) = \pi(x) - \frac{\vartheta(x)}{\log x}. \qquad (23)$$

This completes the verification. ∎

The corresponding representations over arithmetic progressions are as follow.

**Theorem 17.** Let $x > 2$, let $a, q$ be relatively prime numbers, and let $\chi$ be a nontrivial character modulo $q$. Then

(i) $\pi(x,q,a) = \dfrac{\vartheta(x,q,a)}{\log x} + \int_2^x \dfrac{\vartheta(t,q,a)dt}{t \ln^2 t}$, (ii) $\vartheta(x,q,a) = \pi(x,q,a)\log x - \int_2^x \dfrac{\pi(t,q,a)dt}{t \ln^2 t}$, (24)

(iii) $\displaystyle\sum_{p \leq x} \dfrac{\chi(p)}{p} = \dfrac{\vartheta(x,\chi)}{x \log x} + \int_2^x \vartheta(t,\chi)\dfrac{1+\ln t}{t^2 \ln^2 t}dt$, (iv) $\displaystyle\sum_{p \leq x,\, p \equiv a \bmod q} 1/p = \dfrac{\pi(x,q,a)}{x} + \int_2^x \dfrac{\pi(t,q,a)}{t^2}dt$.

These representations provide effective conversions from $\vartheta(x), \psi(x)$ to $\pi(x)$ and vice versa.

## 3. Main Results

The circle of equivalent concepts error term of the prime counting function, primes in short intervals, the distribution of the zeros of the zeta function on the critical strip, and so on, is explored here. It will be shown that advances in the theory of primes in short intervals lead to equivalent advances on the error term of the primes counting function and related concepts.

**Theorem 1.** Let $x > 1$ be a sufficiently large number. Then

$$\pi(x) = li(x) + O(x^{7/12-\varepsilon(x)}),\qquad (25)$$

where $0 \leq \varepsilon(x) < 1/12$ and $\varepsilon(x) \to 0$ as $x \to \infty$.

Proof: Using $\vartheta(x) = x + ax^{7/12-\varepsilon(x)}$, see Theorem 13, and applying the theta to pi conversion formula, Theorem 16-i, return

$$\begin{aligned}\pi(x) &= \dfrac{\vartheta(x)}{\log x} + \int_2^x \dfrac{\vartheta(t)}{t \log^2 t}dt \\ &= \dfrac{x + ax^{7/12-\varepsilon(x)}}{\log x} + \int_2^x \dfrac{t + at^{7/12-\varepsilon(x)}}{t \log^2 t}dt \\ &= \dfrac{x}{\log x} + \int_2^x \dfrac{1}{\log^2 t}dt + \dfrac{ax^{7/12-\varepsilon(x)}}{\log x} + a\int_2^x \dfrac{1}{t^{5/12+\varepsilon(x)} \log^2 t}dt \\ &= li(x) + \dfrac{ax^{7/12-\varepsilon(x)}}{\log x} + a\int_2^x \dfrac{1}{t^{5/12+\varepsilon(x)} \log^2 t}dt \\ &= li(x) + O(x^{7/12-\varepsilon(x)}),\end{aligned} \qquad (26)$$

where $0 \leq \varepsilon(x) < 1/12$, and $\varepsilon(x) \to 0$ as $x \to \infty$. Here the error term $E(x) = O(x^{7/12-\varepsilon(x)})$ absorbs both the middle term and the integral in the fourth line, which is of order $O(x^{7/12-\varepsilon(x)}/\log^2 x)$. ∎

## 4. Applications

Last but not least, other applications of Theorem 13 to the distribution of the zeros of the zeta function on the critical strip and Mertens function are considered here.

A *zero-free region* of the zeta function is a complex half plane of the form $\mathcal{H} = \{ s = \sigma + it \in \mathbb{C} : \zeta(s) \neq 0 \text{ for } \mathfrak{Re}(s) = $ $\sigma > 1 - \sigma_0(t)$, $t \in \mathbb{R}\}$, where $\sigma_0(t) > 0$ is a real valued function or some other definition. The study of the zero-free regions of the zeta function involves elaborate estimates of $\zeta(s)$, $\zeta'(s)$, and the exponential sum $\sum_{x<t<2x} n^{it}$, see [ES, p.88], [RM, p. 94], [IV, p. 143], and [FD] etc.

A different approach explored in this section revolves around a joint application of theta function and estimates of primes in short intervals.

**Theorem 18.** Let $s = \sigma + it \in \mathbb{C}$ be a complex number, and define the complex half plane $\mathcal{H}_1 = \{ s \in \mathbb{C} : \Re(s) > 7/12 \}$. Then $\zeta(s) \neq 0$ for all $s \notin \mathfrak{I}_0$ such that $\Re(s) = \sigma > 0$. In particular, the nontrivial zeros are confined to the subcritical strip $\mathcal{E}_1 = \{ s \in \mathbb{C} : 5/12 \leq \Re(s) \leq 7/12 \}$.

Proof: First, the claim is proved for $\sigma = \Re(s) > 7/12$, and then by the functional equation $\zeta(s) = \zeta(1 - s)$, it follows that $\zeta(s) \neq 0$ for $\sigma = \Re(s) < 5/12$. As given in [IV, p. 305], the logarithm derivative of the zeta function has the integral representation

$$-\frac{\zeta'(s)}{\zeta(s)} = \sum_{n=1}^{\infty} \Lambda(n) n^{-s} = s \int_1^{\infty} \frac{\psi(t)}{t^{s+1}} dt, \qquad (27)$$

for $\sigma > 1$. Substituting $\psi(x) = x + cx^{7/12 - \varepsilon(x)}$, $c > 0$ constant, returns

$$-\frac{\zeta'(s)}{\zeta(s)} = \frac{s}{s-1} + cs \int_1^{\infty} \frac{t^{7/12 - \varepsilon(t)}}{t^{s+1}} dt \qquad (28)$$

Clearly, for $\sigma > 7/12$, the integral is a holomorphic function of the complex number $s = \sigma + i\gamma$. Consequently, the zeta function cannot vanish on the complex half plane $\{ s = \sigma + i\gamma \in \mathbb{C} : \sigma > 7/12 \}$. ∎

Furthermore, since the zeros of the zeta function are fixed and independent of the value of $x > 1$, using a result in [SG], see Theorem 4, it is likely that the zeta function has all its nontrivial zeros on the critical strip subcritical strip $\mathcal{E}_\delta = \{ s \in \mathbb{C} : 1/2 - \delta \leq \Re(s) \leq 1/2 + \delta \}$, $\delta > 0$. The argument is the same as the one below.

**Theorem 19.** Let $s = \sigma + it \in \mathbb{C}$ be a complex number, and define the complex half plane $\mathcal{H}_1 = \{ s \in \mathbb{C} : \Re(s) > 21/40 \}$. Then $\zeta(s) \neq 0$ for all $s \in \mathcal{H}_1$. In particular, the nontrivial zeros of the zeta function are confined to the subcritical strip $\mathcal{E}_1 = \{ s \in \mathbb{C} : 19/40 \leq \Re(s) \leq 21/40 \}$.

Proof: First, the claim is proved for $\sigma = \Re(s) > 21/40$, and then by the functional equation $\zeta(s) = \zeta(1 - s)$, it follows that $\zeta(s) \neq 0$ for all complex numbers such that $0 < \sigma < 19/40$ too. For $x \geq x_0$, the theta difference over the short interval $(x, x + y]$, $y = x^\beta$, $\beta > 0$, is given by

$$\begin{aligned}\vartheta(x + y) - \vartheta(x) &= x + y + O((x + y)^\theta \log^2(x + y)) - \left(x + O(x^\theta \log^2 x)\right) \\ &= y + O(x^\theta \log^2 x),\end{aligned} \qquad (29)$$

where the implied constants are distinct, see Theorem 11. By Theorem 2, the interval $(x, x + x^\beta]$ contains primes whenever $\beta \geq 21/40$, this fact immediately yields

$$\log x < \sum_{x<p\leq x+y} \log p \leq \frac{2y}{\log(y)+3.53}\log(x+y) \leq c_1 y, \qquad (30)$$

the middle inequality follows from Theorem 3, the right inequality follows from $y = x^\beta$, $\beta > 0$, and some constant $c_1 > 0$. Next use the identity $\vartheta(x+y) - \vartheta(x) = \sum_{x<p\leq x+y} \log p$ to combine (1) and (2) into

$$\log x < y + O(x^\theta \log^2 x) \leq c_1 y. \qquad (31)$$

Dividing by $y \geq x^\beta = x^{21/40}$ across the board returns

$$0 \leq x^{-\beta} \log x < 1 + O(x^{\theta-\beta} \log^2 x) \leq c_3. \qquad (32)$$

Clearly, since both sides of these inequalities are bounded by nonnegative constants for all $x \geq x_0$, it implies that $.525 = \beta > \theta$. ∎

In summary, the nontrivial zeros of the analytic continuation $\zeta(s) = (1-2^{1-s})^{-1} \sum_{n\geq 1} (-1)^{n-1} n^{-s}$, $\Re(s) > 0$, of the zeta function $\zeta(s) = \sum_{n\geq 1} n^{-s}$, $\Re(s) > 1$, are on the subcritical strip $\mathcal{E}_1 = \{ s \in \mathbb{C} : 19/40 \leq \Re(s) \leq 21/40 \}$. Furthermore, since the zeros of the zeta function are fixed and independent of the value of $x > 1$, a known result on the existence and density of prime numbers in almost all short intervals, see Theorem 4, probably implies that the zeta function has all its nontrivial zeros confined to the subcritical strip $\mathcal{E}_\delta = \{ s \in \mathbb{C} : 1/2 - \delta \leq \Re(s) \leq 1/2 + \delta \}$, $\delta > 0$.

By means of the explicit formula, the zero-free region $\{ \sigma > 1 - c(\log t)^{-\beta} \}$ is mapped to an approximation of the theta function of the shape $\vartheta(x) = x + O(xe^{-c(\log x)^{1/(1+\beta)}})$, see [TM, p.56], [IG, p. 60]. In contrast, the zero-free region $\{ \sigma > \sigma_0 \}$ has a simpler correspondence, namely, $\vartheta(x) = x + O(x^{\sigma_0})$.

**Corollary 20.** Let $x \geq x_0$. Then
(i) $\psi(x) = x + O(x^{21/40})$.      (ii) $\vartheta(x) = x + O(x^{21/40})$.

**Mertens Function.** The Mertens function $M(x) = \sum_{n\leq x} \mu(n)$ is known to satisfy the inequality $|M(x)| > x^{1/2}$. Specifically, it was shown that

$$\liminf_{x\to\infty} \frac{M(x)}{\sqrt{x}} < -1.009 \quad \text{and} \quad \limsup_{x\to\infty} \frac{M(x)}{\sqrt{x}} > 1.06, \qquad (33)$$

see [OT], and it is believed to satisfy $|M(x)| \leq cx^{1/2}$, for some constant $c > 1$. The current best estimate of this function has the shape stated below.

**Theorem 21.** For real numbers $x \geq x_0$, the Mertens function satisfies $\left|\sum_{n\leq x} \mu(n)\right| = O(xe^{-c(\log x)^{1/2}})$, where $c > 0$ is an absolute constant.

An estimate of this function does not have a main term since its generating function is holomorphic.

**Theorem 22.** For all sufficiently large numbers $x > 1$, the Mertens function satisfies $\left|\sum_{n\leq x} \mu(n)\right| \leq x^{7/12}$.

*Proof*: The inverse $1/\zeta(s)$ of the zeta function $\zeta(s)$ for $\sigma > 1$ is derived from the product of power series

$$\zeta(s)\sum_{n=1}^{\infty}\frac{\mu(n)}{n^s} = \sum_{n=1}^{\infty}\frac{1}{n^s}\sum_{n=1}^{\infty}\frac{\mu(n)}{n^s} = \sum_{n=1}^{\infty}\left(\sum_{d|n}\mu(d)\right)n^{-s} = 1, \tag{34}$$

where $s = \sigma + it \in \mathbb{C}$. Moreover, the integral representation of $1/\zeta(s)$ is given by

$$\begin{aligned}\frac{1}{\zeta(s)} &= \sum_{n=1}^{\infty}\frac{\mu(n)}{n^s} = \sum_{n=1}^{\infty}\frac{M(n)-M(n-1)}{n^s} = \sum_{n=1}^{\infty}M(n)\left(\frac{1}{n^s}-\frac{1}{(n-1)^s}\right) \\ &= \sum_{n=1}^{\infty}M(n)\left(s\int_n^{n+1}\frac{dx}{x^{s+1}}\right) = s\int_1^{\infty}\frac{M(x)}{x^{s+1}}dx,\end{aligned} \tag{35}$$

where the third equality follows from a change of index $n \to n+1$ in the second fraction only. Moreover, since the zeta function $\zeta(s)$ does not vanish on the complex half plane $\mathcal{H}_1 = \{ s = \sigma + it \in \mathbb{C} : \sigma > 7/12 \}$, see the Theorem 18, it implies that $1/\zeta(s)$ is a bounded function of the complex number $s = \sigma + it \in \mathfrak{I}$. This in turn implies that the integral is a holomorphic function of the complex variable $s = \sigma + it \in \mathbb{C}$, $\sigma > 7/12$. Therefore, it follows that it absolute value satisfies $|M(x)| \leq x^{7/12}$. ∎

The quasirandom nature of the function $M(x) = \sum_{n \leq x}\mu(n)$, where $\mu(n) = -1, 0, 1$, in tandem with the oscillatory term $x^{s+1} = e^{-(\sigma+it)\log x}$ can cause significant cancellation on the integral. Due to these possible cancellations in the integral

$$\frac{1}{\zeta(s)} = s\int_1^{\infty}\frac{M(x)}{x^{s+1}}dx = s\int_1^{\infty}M(x)e^{-(\sigma+it)\log x}dx, \tag{36}$$

it is likely that $|M(x)| \leq cx^{1/2}$ for some constant $c > 1.06$. For a truly random variable $S(x) = \sum_{n \leq x}x(n)$, where $x(n) = -1, 1$ is an independent random variable, the limit superior is known to satisfies the (Kolmogorov) expression

$$\limsup_{x \to \infty}\frac{S(x)}{\sqrt{.5x\log\log x}} = 1, \tag{37}$$

with probability 1. Amazingly, the function $Q(x) = \sum_{n \leq x}|\mu(n)| = 6\pi^{-2}x + O(x^{1/2})$, which counts the number of squarefree integers $n \leq x$, is much easier to determine than $M(x)$.

## 5. Appendix
**Power Series**

A) $(1+z)^\tau = 1 + z + \tau\frac{z^2}{2!} + \tau(\tau-1)\frac{z^3}{3!} + \cdots = 1 + z + O(z^2)$, for $|z|<1$ and $\tau > 0$. (38)

B) $\log(1+z) = z + \frac{z^2}{2} + \frac{z^3}{3} + \cdots = z + O(z^2)$, for $|z|<1$. (39)

C) $(x+y)^N = x^N(1+y/x)^N = x^N(1+Ny/x+O(y^2/x^2))$, for a fixed $N > 0$. (40)

D) $\log^N(x+y) = \log^N x + N(\log^{N-1} x)\log(1+y/x) + \cdots + \log^N(1+y/x) = \log^N x + O((\log^{N-1} x)\log(1+y/x))$, for a fixed $N > 0$. (41)

**Integrals**

The real and complex logarithm integrals are defined by

$$li(x) = \int_2^x \frac{dt}{\ln t} \text{ for } x \in \mathbb{R}, \text{ and } li(e^z) = \int_{-\infty+iv}^z e^t\frac{dt}{t} \text{ for } z = u+iv \in \mathbb{C}, v \neq 0$$ (42)

respectively. The complex logarithm integral appears in the Riemann explicit formula

$$\pi(x) = li(x) - \sum_\rho \left(li(x^\rho) + li(x^{1-\rho})\right) + \log\xi(0) + \int_x^\infty \frac{dt}{t(t^2-1)\log t},$$ (43)

which clearly shows the influence of the zeros of the zeta function on the prime counting function $\pi(x)$, see also (5).

B) $\int_2^x \frac{dx}{x^\nu \log^C x} = \frac{x^{1-\nu}}{\log^C x} + \int_2^x \left(\nu + \frac{C}{\log x}\right)\frac{dx}{x^\nu \log^C x} = O\left(\frac{x^{1-\nu}}{\log^C x}\right)$ for any constants $C, \nu > 0$ and $x > 2$. (44)

**Proofs**

*Case* 1. The constants $c_0 \neq c_1$, and $\vartheta(x) = x + c_0 x/f(x)$, where $f(x) = (\log x)^B$, $B > 0$ is a constant. The theta difference over the short interval $(x, x+y]$ is given by

$$\vartheta(x+y) - \vartheta(x) = y + \frac{c_1(x+y)}{\log^B(x+y)} - \left(x - \frac{c_0 x}{\log^B x}\right)$$
$$= y\left(1 + \frac{c_1}{\log^B(x+y)}\right) + x\left(\frac{c_1}{\log^B(x+y)} - \frac{c_0}{\log^B x}\right)$$ (45)

Use the power series expansion (D) to simplify the difference

$$\frac{c_1}{\log^B(x+y)} - \frac{c_0}{\log^B x} = \frac{c_1\log^B x - c_0\log^B(x+y)}{\log^B x \log^B(x+y)} = \frac{c_1\log^B x - c_0(\log^B x + O(\log^{B-1} x \log(1+y/x)))}{\log^B x \log^B(x+y)}$$
$$= \frac{c_1 - c_0 + O(\log(1+y/x))}{\log^{B+1}(x+y)}$$ (46)

Replacing it back yields

$$\vartheta(x+y) - \vartheta(x) = y\left(1 + \frac{c_1}{\log^B(x+y)}\right) + x\left(\frac{c_1 - c_0 + O(\log(1+y/x))}{\log^{B+1}(x+y)}\right). \tag{47}$$

Next rewrite it as $\vartheta(x+y) - \vartheta(x) = \sum_{x<p\leq x+y} \log p$ and apply Theorem 3 to estimate the inequalities

$$\frac{c_2 y}{\log(x+y)} \log x \leq \sum_{x<p\leq x+y} \log p \leq \frac{c_3 y}{\log(x+y)} \log(x+y). \tag{48}$$

Upon division by $y > c_4 x^{7/12 - \varepsilon(x)}$, the penultimate equation becomes

$$c_2\left(1 - \frac{y}{x \log x}\right) \leq 1 + \frac{c_1}{\log^B(x+y)} + x^{5/12 + \varepsilon(x)}\left(\frac{c_1 - c_0 + O(\log(1+y/x))}{\log^{B+1}(x+y)}\right) \leq c_3. \tag{49}$$

But since both the left side and the right side of these inequalities are constants for sufficiently large $x > 0$, this is a contraction for any $B > 0$. Thus, it quickly follows that $f(x) \neq (\log x)^B$.